\documentstyle[12pt]{article}
\usepackage {latexsym, amsfonts, emlines, emlines2}
\begin{document}
\begin{center} {\small\it Algebra and Logic, Vol. 38, No. 4,
1999, 259--276.} \hfill \vspace*{2\baselineskip}

{\large\bf A QUANTUM ANALOG OF THE\\ 
POINCARE--BIRKHOFF--WITT THEOREM}
\hspace*{\fill} \\ [2mm] \hspace*{6mm} 

{\large\bf V. K.
Kharchenko}
\footnote{Supported by the National Society of Researchers,
M\'exico (SNI, exp. 18740, 1997-2000).}

\

Translated from {
Algebra i Logika, Vol. 38, No. 4, pp. 476-507, July-August, 1999.
Original article submitted June~29, 1998.} \end{center}

\begin{quote}\it{We reduce the basis construction problem for
Hopf algebras generated by skew-primitive semi-invariants
to a study of special elements, called
``super-letters,'' which are defined by Shirshov standard words. 
In this way we show that above Hopf algebras always have sets of
PBW-generators (``hard'' super-letters).
It is shown also that these Hopf algebras having not more than finitely
many ``hard'' super-letters share some of the properties of
universal enveloping algebras of finite-dimensional Lie algebras. The
background for the proofs is the construction of a filtration such
that the associated graded algebra is obtained by iterating the skew
polynomials construction, possibly followed with
factorization.}\end{quote}

\begin{center}
{\bf INTRODUCTION}
\end{center}
\smallskip

In this article we deal with the basis construction problem for
character Hopf algebras, i.e., for the Hopf algebras generated by
skew primitive semi-invariants and by an Abelian group of all
group-like elements. These algebras constitute an important class
actively studied within the frames of the quantum group theory. The
class includes all known to the date quantizations with the coalgebra
structure of a Lie algebra, and probably we may think of it as an
abstractly defined class of all ``quantum'' universal enveloping
algebras. In line with this approach, ``quantum'' Lie algebras are
couched in terms of spaces of all skew primitive elements of the
character Hopf algebras endowed with the natural structure of an
Yetter--Drinfeld module and equipped with partial quantum
operations; see [1].

In the present article the basis construction problem will be
reduced to treating special elements defined by Shirshov standard
words, which we call ``super-letters.'' The main result, Theorem 2,
states that the set of all monotonic restricted words in ``hard''
super-letters constitute a basis for a Hopf algebra. If the Hopf
algebra is generated by ordinary primitive elements, the set of all
``hard'' super-letters constitute a basis for the Lie algebra of all
primitive elements. By this token, Theorem 2 may be conceived of as
one of the possible quantum analogs for the Poincare--Birkhoff--Witt
theorem.

The proof and the statement of the main theorem are based on
Shirshov's combinatorial method, originally developed for Lie
algebras. Using this method, Shirshov solved a number of important
problems in combinatorial theory of Lie algebras. Among them are the
characterization problem for a free Lie algebra over an arbitrary
operator ring [3], the basis construction problem for a free Lie
algebra over a field [5] (which is an outgrowth of M. Hall's ideas
in [4]), and the equality problem for Lie algebras with one defining
relation or with finitely many homogeneous relations [6].
Independently, most fundamental concepts of that method were
pronounced in [7] where the basis construction problem was dealt
with for dual groups of the lower central series of a finitely
generated free group.

A weak point in our modification of Shirshov's method is that
essential use will be made of the so-called ``through'' ordering of
words, standard words, and super-words, for which the set of all
standard words (super-letters) may be not completely ordered. For
this reason, the main theorem is proved only for finitely generated
Hopf algebras. In this connection, it is worth mentioning that
Shirshov's original method does not presume the use of a ``through"
ordering only. What it calls for is a weak restriction on the order:
the end of a standard word should be less than the word itself. An
example is M. Hall's ordering from [4], which is in fact also used
in the present article. However, we opt to not bring in both of the
orders, to avoid (or at least minimize) misunderstanding. The reason
why we do not use  the M. Hall's ordering as the  main to our
reasoning is because Lemma 8 becomes almost uninformative in this
case.

In order to generalize the main theorem to the case of infinitely
generated algebras, instead of searching suitable orders, one may
apply the famous local method by Mal'tsev (cf. [8]), whereby the
proof of the theorem reduces to a logical analysis of its
formulation.

The main theorem can also be used to construct bases for known
quantizations of Lie algebras. For the Drinfeld--Jimbo
quantizations, such were constructed in Rosso [9], Yamane [10],
Lusztig [11, 12], and Kashiwara [13]. Curiously, no one of the
methods by Rosso, Yamane, Lusztig, or Kashiwara presupposes that
active use be made of the coalgebraic structure --- instead --- they
all presume a detailed treatment of the algebraic. At the
same time, the coproduct in Drinfeld--Jimbo quantizations, and also
in every pointed Hopf algebra (cf. [14]), consists only of a ``skew
primitive'' leading part and a linear tensor combination of a lesser
degree. --- This opens up unbounded prospects for inductive proofs.

An approach attempted here is aimed at a study of effects brought
about by the existence of a coproduct. The Poincare--Birkhoff--Witt
theorem (PBW-theorem) can also be proved in terms of a coproduct,
provided that a given Lie algebra is presupposed to be embedded in a
(cocommutative connected) Hopf algebra. This was in fact done in
Milnor and Moore [15, Secs. 5 and 6]. The PBW-theorem in the
Milnor-Moore form carries no information about primitive elements
(the given Lie algebra) but gives a complete solution to the basis
construction problem for a Hopf algebra modulo its solution for a
given Lie algebra. The mentioned above decreasing process 
(unlike detailed algebraic accounts) has sharply 
delineated boundaries of application --- it
cannot give any information about the structure of a set of skew
primitive elements (that is  of the structure of the quantum Lie
algebra itself). Therefore, it might be interesting to investigate these
sets in known quantizations  as quantum Lie algebras, that is together
with all partial quantum operations over them [1].

In Sec. 1, we introduce basic notions and give a formulation of
Shirshov's theorem [3, Lemma 1] needed for our further
constructions. All statements under this section were proved by
Shirshov  sometimes in a more general
form. In a slightly different guise, some of them were discovered 
independently of Shirshov in [7].

In Sec. 2, we replace the classical commutator with a skew
commutator whose ``curvature'' depends on the parameters of
specified elements in approximately the same way as it does in color
Lie super-algebras. In our case, however, the bicharacter is not
assumed symmetric. Still, identity (8), which is analogous to the
Jacobi identity, is valid. And so are derivative identities (9) and
(11), which link the skew commutator and the basic product. The bulk
of the information needed is given in Lemmas 6 and 8, in which two
decreasing processes are described. One is an analog of the
Hall--Shirshov construction for nonassociative words and the other
is concerned with a coproduct in the way mentioned above.

In Sec. 3, we pass from quantum variables to arbitrary skew
primitive generators, and using the two above-mentioned decreasing
processes, prove the main result, Theorem 2. On this theorem, each
character Hopf algebra has the same basis as the universal
enveloping algebra of a (restricted) Lie algebra. The role of a
basis for the Lie algebra is played by special elements, which we
call hard super-letters. The main lemmas are stated in such a way as
to fit in dealing with skew primitive elements.

In Sec. 4, we derive some immediate consequences of the main
theorem. In particular, it is shown that character Hopf algebras
having not more than finitely many hard super-letters share some of
the properties of universal enveloping algebras of finite-dimensional Lie
algebras. The background for our proofs is the construction of a
filtration such that the associated graded algebra is obtained by
iterating the skew polynomials construction, possibly followed with
factorization. Note also that the main theorem, as well as its
corollaries, remain true for $(G,\lambda)$-graded Hopf algebras
and for braided bigraded Hopf algebras. In this event a group $G$
merely defines a grading, but the algebra in question does not
itself contain the $G$. Therefore, additional restrictions on a
group are unnecessary.

Finally, the quantum Serre relations can be expressed in terms of
some super-letters being equal to zero. If, in these super-letters,
we replace the skew commutator operation with the classical one then
the original Serre relations will appear; see [1, Thm. 6.1].
Therefore, it seems absolutely realistic that all hard super-letters
of the Drinfeld--Jimbo quantized enveloping algebras arise from a
suitable basis for the Lie algebra by merely replacing the commutator
with the  skew commutator. This is likely to be true not only for the
case of Drinfeld--Jimbo quantizations.

\medskip
\begin{center}
{\bf 1. SHIRSHOV STANDARD WORDS}
\end{center}
\smallskip

Let $x_1,\ldots ,x_n$ be a set of variables. Consider this set as
an alphabet. On a set of all words in this alphabet, define the
lexicographical order such that $x_1>x_2>\ldots >x_n$. This means
that two words $v$ and $w$ are compared by moving from left to
right until the first distinct letter is encountered. If not,
i.e., one of the words is the beginning of the other, then a shorter
word is assumed to be greater than the longer (as is common practice
in dictionaries). For example, all words of length at most two in
two variables respect the following order:
\begin{equation}
x_1>x_1^2>x_1x_2>x_2>x_2x_1>x_2^2.
\label{ex2}
\end{equation}

\noindent This order is stable under left multiplication and
unstable under right. Nevertheless, if $u>v$ and $u$ is not the
beginning of $v$, then the inequality is preserved under right
multiplication, even by different words: $uw>vt$.

Every noncommutative polynomial $f$ in $x_1,\ldots,x_n$ is a
linear combination of words $f=\sum \alpha _iu_i$. By
$\overline{f}$ we denote a {\it leading word} which occurs in this
decomposition with a nonzero coefficient. In the general case the
leading word of a product does not equal the product of leading
words of the factors. For example, if $f=x_1+x_1x_2$ and
$g=x_3+x_3x_2$ then $\overline{fg}=x_1x_2x_3\neq
x_1x_3=\overline{f}\overline{g}$. But if the leading word of $f$
is not the beginning of any other word in $f$, then
\begin{equation}
\overline{fg}=\overline{f}\overline{g}.
\label{sta}
\end{equation}

\noindent Indeed, the inequalities $\overline{f}>u_j$ can be
multiplied from the right by (possibly distinct) elements
$\overline{f}v_k>u_jv_s$. In particular, if $f$ is an
homogeneous polynomial, i.e., all words $u_i$ have the same
length, then formula (2) is true.

The set of all words is not completely ordered since there exist
infinite decreasing chains --- for instance,
\begin{equation}
x_1>x_1^2>x_1^3>\ldots >x_1^n>\ldots \,.
\label{ex1}
\end{equation}

\noindent Yet, all of its finite subsets are completely ordered.
This will allow us to use induction on the leading word, provided
that bounds are set on the lengths of words, $l(v)$, or on degrees
of the polynomials under consideration.

\smallskip
{\bf Definition 1.} A word $u$ is called {\it standard} (in the
sense of Shirshov) if, for each representation $u=u_1u_2$, where
$u_1$ and $u_2$ are nonempty words, the inequality $u>u_2u_1$
holds. For example, in (3), there is only one standard word,
namely, $x_1$, and in (1), there are three: $x_1$, $x_1x_2$,
and $x_2$.

\smallskip
{\bf LEMMA 1.} {\it If $u=sv$ is a standard word and $s$ is nonempty
then $v$ is not the beginning of} $u$.

{\bf Proof.} Suppose, to the contrary, that $u=vs^{\prime }$. By
the definition of a standard word, we then have $sv=vs^{\prime
}>s^{\prime }v$, i.e., $s>s^{\prime }$. Similarly, $vs^{\prime
}=sv>vs$, whence $s^{\prime }>s$, a contradiction.

\smallskip
{\bf LEMMA 2.} {\it A word $u$ is standard if and only if it is greater
than any one of its endings}.

{\bf Proof.} If the word $u$ is standard and $u=vv_1$ then
$vv_1>v_1v$ by definition. By the previous lemma, $v_1$ is not
the beginning of $vv_1$, and hence $u=vv_1$ and $v_1v$ differ
already in their first $l(v_1)$ letters. Therefore, $u>v_1$.
Conversely, if $u=u_1u_2$ and $u>u_2$ then $u$ is not the
beginning of $u_2$, and so the above inequality remains true
under right multiplication of the right hand side by $u_1$.

\smallskip
{\bf LEMMA 3.} {\it If $u$ and $v$, $u>v$, 
are standard words  then} $u^h>v$.

{\bf Proof.} If $u$ is not the beginning of $v$ then $u>v$ can be
multiplied from the right by different words. Suppose that
$v=u^kv^{\prime }$ and that $v^{\prime }$ does not begin with $u$.
If $k\geq h$ then $u^h>v$ as a beginning. If 
$k<h$ then $v^{\prime }$ is nonempty and
$v^{\prime }<v<u$. It follows that $v=u^k\cdot
v^{\prime }<u^k\cdot u\cdot u^{h-k-1}=u^h$.

\smallskip
{\bf LEMMA 4.} {\it Let $u$ and $u_1$ be standard words such that
$u=u_3u_2$ and $u_2>u_1$. Then}
\begin{equation}
uu_1>u_3u_1,\quad uu_1>u_2u_1.
\label{for1}
\end{equation}

{\bf Proof.} First we show that $u_2u_1>u_1$. If $u_1$ does not
begin with $u_2$, the inequality follows immediately from
$u_2>u_1$. Assume that $u_1=u_2^k\cdot u_1^{\prime }$ and $u_2$
is not the beginning of $u_1^{\prime }$. Since $u_1$ is standard,
we have $u_2^ku_1^{\prime }>u_2^{k-1}u_1^{\prime }$, i.e.,
$u_2u_1^{\prime }>u_1^{\prime }$. Hence $u_2u_1=u_2^k\cdot
u_2u_1^{\prime }> u_2^k\cdot u_1^{\prime }=u_1$. Multiplying this
inequality from the left by $u_3$ yields the first inequality
required. Consider the second. Since $u$ is a standard word,
$u_3u_2>u_2$ by Lemma 2, and $u_3u_2$ is of course not the beginning
of $u_2$. We can therefore multiply the latter inequality from the
right by $u_1$.

Recall that a {\it nonassociative word} is one where $[\, , \,]$ are
somehow arranged to show how multiplication applies. The
set of nonassociative words can be defined inductively by the
following axioms:

(1) all letters are nonassociative words;

(2) if $[v]$ and $[w]$ are nonassociative words then $\bigl[
[v][w]\bigr] $ is a nonassociative word;

(3) there are no other nonassociative words.

\smallskip
{\bf Definition 2.} A nonassociative word $[u]$ is said to be {\it
standard} (in the sense of Shirshov) if:

(1) an (associative) word $u$ obtained from this word by
removing the brackets is standard;

(2) if $[u]=\bigl[ [v][w]\bigr]$ then $[v]$ and $[w]$ are
standard nonassociative words;

(3) if $[u]=\bigl[ \bigl[ [v_1][v_2]\bigr] [w]\bigr] $ then
$v_2\leq w$.

\smallskip
{\bf The Shirshov theorem} (cf. [3, Lemma 1]). Each standard word
can be uniquely bracketed so that the resulting nonassociative word
is standard.

This theorem, combined with the inductive definition of a set of all
nonassociative words, immediately implies that every standard
(associative) word $u$ has a decomposition $u=vw$, where $v>w$
and $v$ and $w$ are standard. Yet, for the associative
decomposition (as distinct from nonassociative one), the words $v$ and
$w$ are not defined uniquely. The factors $v$ and $w$ in the
nonassociative decomposition $[u]=[[v][w]]$, we note, can be
defined to be standard words such that $u=vw$, where $v$ has a
least possible length; see [16].

\medskip
\begin{center}
{\bf 2. DECOMPOSITION OF QUANTUM POLYNOMIALS  INTO LINEAR COMBINATIONS OF 
MONOTONIC SUPER-WORDS }
\end{center}
\smallskip

Let $x_1,\ldots,x_n$ be {\it quantum variables}, i.e., associated
with each letter $x_i$ are an element $g_i$ of a fixed Abelian
group $G$ and a character $\chi ^i:G\rightarrow {\bf k}^*$. For
every word $u$, denote by $g_u$ an element of the group $G$
which results from $u$ by replacing each occurrence of the letter
$x_i$ with $g_i$. This group-like element is denoted also by
$G(u)$, provided that $u$ is an unwieldy expression. Likewise,
by $\chi ^u$ we denote a character which results from $u$ by
replacing all $x_i$ with $\chi ^i$. For a pair of words $u$ and
$v$, put
\begin{equation}
p_{uv}=\chi ^u(g_v).
\label{p1}
\end{equation}

\noindent Obviously, the following equalities hold:
\begin{equation}
p_{uu_1\ v}=p_{uv}p_{u_1v},\quad p_{u\ vv_1}=p_{uv}p_{uv_1},
\label{p2}
\end{equation}

\noindent that is the operator $p$ is a bicharacter defined on a
semigroup of all words. Sometimes we denote this operator by
$p(u,v)$. Define a bilinear termal operation, a skew commutator,
on a set of all quantum polynomials by setting
\begin{equation}
[u,v]=uv-p_{uv}vu.
\label{op1}
\end{equation}

\noindent This satisfies the following identity:
\begin{equation}
[[u,v],w]=[u,[v,w]]+p_{wv}^{-1}[[u,w],v]+(p_{vw}-p_{wv}^{-1})[u,w]\cdot v,
\label{Jak1}
\end{equation}

\noindent which is similar to the Jacobi identity, where $cdot$
stands for usual multiplication in a free algebra, and which
can be easily verified by direct computations using (7). Likewise,
the following formulas of skew derivations, by which the skew
commutator is linked to multiplication, are valid:
\begin{equation}
[u,v\cdot w]=[uv]\cdot w+p_{uv}v\cdot [uw],
\label{dif1}
\end{equation}
\begin{equation}
[u\cdot v,w]=p_{vw}[uw]\cdot v+u\cdot [vw].
\label{dif2}
\end{equation}

\smallskip
{\bf Definition 3.} A {\it super-letter} is a polynomial equal to a
standard nonassociative word with brackets defined as in operation
(7).

By the Shirshov theorem, every standard word $u$ is associated with
a super-letter $[u]$. If we remove the brackets in $[u]$ as is
done in definition (7) we obtain an homogeneous polynomial whose
leading word is equal to $u,$ and this leading word occurs in the decomposition
of $[u]$ with coefficient 1. This is easily verified by induction
on the degree. Indeed, if $[u]=[[v][w]]$ then the super-letter
$[u]$ is equal to $[v][w]-p_{uw}[w][v]$. By the induction
hypothesis, $[v]$ and $[w]$ are homogeneous polynomials with the
leading words $v$ and $w$, respectively. Therefore, the leading
word of the first summand equals $vw$ and has coefficient 1; the
leading word of the second equals $wv$ and is less than $vw$ by
definition.

Thus, in correspondence with distinct standard words $u$ and $v$
are distinct super-letters $[u]$ and $[v]$, and the order on a
set of super-letters can be defined as follows:
\begin{equation}
[u]>[v]\iff u>v.
\label{por1}
\end{equation}

\smallskip
{\bf Definition 4.} A word in super-letters is called a {\it super-word}. A super-word is said to be {\it
monotonic} if it has the form \begin{equation}
W=[u_1]^{k_1}[u_2]^{k_2}\cdots [u_m]^{k_m},
\label{mon}
\end{equation}

\noindent where $u_1<u_2<\ldots <u_m$.

We recall that the {\it constitution} of $u$ is a sequence of integers
$(m_1,m_2,\ldots, m_n)$ such that $u$ has degree $m_1$ in $x_1$,
degree $m_2$ in $x_2$, etc. Since super-letters and super-words
are homogeneous in each of the variables, their constitutions can be 
defined in the obvious manner. Because $G$ is commutative, the
elements $g_u$ and the characters $\chi ^u$ are the same for all
words of a same constitution. For super-letters and super-words,
therefore, $G(W)=g_{w}$ and $p(U,V)=p_{uv}$ are defined uniquely.

On the set of all super-words, consider a lexicographic order
defined by the ordering of super-letters in (11).

\smallskip
{\bf LEMMA 5.} {\it A monotonic super-word
$W=[w_1]^{k_1}[w_2]^{k_2}\cdots [w_m]^{k_m}$ is greater than a
monotonic super-word $V=[v_1]^{m_1}[v_2]^{m_2}\cdots [v_k]^{m_k}$
if and only if the word $w=w_1^{k_1}w_2^{k_2}\cdots w_m^{k_m}$ is
greater than the word $v= v_1^{m_1}v_2^{m_2}\cdots v_k^{m_k}$.
Moreover, the leading word of the polynomial $W$, when decomposed
into a sum of monomials, equals $w$ and has coefficient 1}.

{\bf Proof.} Let $W>V$. Then $w_1\geq v_1$ in view of the
ordering of super-letters. If $w_1=v_1$, we can remove one
factor from the left of both $V$ and $W$, and then proceed by induction.
Therefore, we will put $w_1>v_1$. If $w_1$ is not the beginning of
$v_1$, then the latter inequality can be multiplied from the right
by suitable distinct elements, which yields $w>v$, as required.
Let $v_1=(w_1^{k_1}w_2^{k_2}\cdots w_{s-1}^{k_{s-1}})w_s^l\cdot
v_1^{\prime }$, where $0\leq l<k_s$. Note that the term
between the parentheses may be missing (in which case $s=1$,
$l>0$), and $w_s$ is not the beginning of $v_1^{\prime }$. If
$v_1^{\prime }$ is a nonempty word, then $v_1^{\prime }<v_1<w_1\leq
w_s$, since $v_1$ is standard. To obtain $v<w$, the inequality
$v_1^{\prime }<w_s$ will be multiplied from the left by one element
$(w_1^{k_1}w_2^{k_2}\cdots w_{s-1}^{k_{s-1}})w_s^l$, and from the
right by (possibly) different elements. If
$v_1^{\prime }$ is the empty word, again we arrive at a
contradiction with $v_1$ being standard. Indeed, if $l>0$, then
the word $v_1$ should be greater than its end $w_s$;
therefore, $w_1>v_1>w_s$, which contradicts the fact that $w_1\leq
w_s$ is valid for all $s\geq 1$. If $l=0$, then $s>1$,
since $v_1$ begins with $w_1$. It follows that $v_1$ is greater
than its end $w_{s-1}$, which is again a contradiction with
$w_1>v_1>w_{s-1}$.

The second part of the lemma follows from the fact that the leading
word of a product of homogeneous polynomials equals the product of
leading words of the factors.

The lemma cannot be extended to the case of nonmonotonic
super-words, for example, $[x_1]\cdot [x_3]>[x_1x_2]$ and
$x_1x_3<x_1x_2$.

\smallskip
{\bf LEMMA 6.} {\it Let $u$ and $u_1$ be standard words and $u>u_1$.
Then the polynomial $[[u][u_1]]$ is a linear combination of
super-words in the super-letters $[w]$ which lie properly between
$[u]$ and $[u_1]$ and are such that $w\leq uu_1$. In this case
the degree of every summand in each of the variables $x_1,\ldots,
x_n$ is equal to a respective degree of} $uu_1$.

{\bf Proof.} If the nonassociative word $[[u][u_1]]$ is standard
then it defines a super-letter $[w]$. In this case $u>w$ since
$u$ is the beginning of $w$, and $w>u_1$ by Lemma 2. In
particular, the lemma is valid if the degrees of $u$ and $u_1$
are equal to 1. And we can therefore proceed by induction on the
length of $uu_1$.

Suppose that our lemma is true if the length of $uu_1$ is less than
$m$. Choose a pair $u,u_1$ with a greatest word $[u]$, so that
the polynomial $[[u][u_1]]$ does not enjoy the required
decomposition and the length of $uu_1$ equals $m$. Then the word
$[[u][u_1]]$ is not standard, i.e., $[u]=[[u_3][u_2]]$ with
$u_2>u_1$. We introduce the notation for super-letters $U_i=[u_i]$,
$i=1,2,3$. By Jacobi identity (8), we can write
\begin{equation}
[[U_3U_2]]U_1]=[U_3[U_2U_1]]+
p_{u_1u_2}^{-1}[[U_3U_1]U_2]+
(p_{u_2u_1}-p_{u_1u_2}^{-1})[U_3U_1]\cdot U_2.
\label{raz}
\end{equation}

\noindent It follows that $u_3>u>u_2>u_1$. By the induction
hypothesis, $[U_3U_1]$ can be represented as $\sum\limits _i \alpha
_i\prod\limits _k[w_{ik}]$, where $u_3>u_3u_1\geq w_{ik}>u_1$.
Using Lemma 4, we obtain $u>uu_1>u_3u_1\geq w_{ik}$, i.e., all
super-letters $[w_{ik}]$ satisfy the requirements of the present
lemma. Furthermore, the word $u$ cannot be the beginning of $u_2$,
and so $u>u_2$ implies $uu_1>u_2$. Thus the super-letter $U_2$,
too, satisfies the requirements. Consequently, the second [in view
of (7)] and third summands of (13) have the required decomposition.

Using the induction hypothesis, for the first summand we obtain
\begin{equation}
[U_2U_1]=\sum _i \beta _i\prod _k[v_{ik}],
\label{fir}
\end{equation}

\noindent where $u_2>u_2u_1\geq v_{ik}>u_1$. By Lemma 4,
$uu_1>u_2u_1\geq v_{ik}$, i.e., the super-letters $[v_{ik}]$
satisfy the conditions of the lemma. Rewrite the first summand
using skew-derivation formula (9), with the first factor replaced by
(14). With this, the first summand turns into a linear
combination of words in the super-letters $[v_{ik}]$ and skew
commutators $[[u_3][v_{ik}]]$. Since $u_3>u>u_2>v_{ik}$ and the
length of $v_{ik}$ does not exceed that of $u_2u_1$, the
induction hypothesis applies to yield
\begin{equation}
[[u_3][v_{ik}]]=\sum _j\gamma _j\prod _t[w_{jt}],
\label{fir1}
\end{equation}

\noindent where $u_3>u_3v_{ik}\geq w_{jt}>v_{ik}$. In this case
$u_2u_1\geq v_{ik}$ implies $uu_1= u_3u_2u_1\geq u_3v_{ik}\geq
w_{jt}$; in addition, $w_{jt}>v_{ik}>u_1$, i.e., the
super-letters $[w_{jt}]$ also satisfy the conditions.

\smallskip
{\bf LEMMA 7.} {\it Every nonmonotonic super-word is a linear combination
of lesser monotonic super-words of a same constitution, whose
super-letters all lie (not strictly) between the greatest and the
least super-letters of a super-word given}.

{\bf Proof.} We proceed by induction on the degree. Whenever
super-letters of a given super-word are rearranged, the degree of a
polynomial remains fixed; therefore, the least super-word of degree
$\leq m$ will be monotonic. Assume that the lemma is true for
super-words of degree $<m$, letting $W=UU_1\cdots U_t$ be a least
super-word of degree $m$ for which our lemma fails. If the
super-word $U_1\cdots U_t$ is not monotonic, by the induction
hypothesis, then, it is a linear combination of lesser monotonic
super-words $W_i$. And we can now apply the induction hypothesis
to $UW_i$. Let
\begin{equation}
W=UU_1^{k_1}\cdots U_t^{k_t},\ \ \ \ U_1<U_2<\ldots <U_t.
\label{na1}
\end{equation}

\noindent If $U \leq U_1$ then $W$ is monotonic, and there is
nothing to prove. Let $U>U_1$. Then
\begin{equation}
W=[UU_1]U_1^{k_1-1}\cdots U_t^{k_t}+
p_{uu_1}U_1UU_1^{k_1-1}\cdots U_t^{k_t}.
\label{na2}
\end{equation}

\noindent The second summand being a super-word is less than
$W$, and so we can write it in the required form. By Lemma 6, the
factor $[UU_1]$ in the first term can be represented as
$\sum\limits _i\alpha _i\prod\limits _j[w_{ij}]$, where the
super-letters $[w_{ij}]$ are less than $U$. Consequently, the
super-letters $\prod\limits _j[w_{ij}]U_1^{k_1-1}\cdots U_t^{k_t}$
are less than $W$, i.e., the first term, and hence also $W$,
will have the required representation.

\smallskip
{\bf THEOREM 1.} {\it The set of all monotonic super-words constitute a
basis for a free algebra} ${\bf k}\langle x_1,\ldots,x_n\rangle $.

{\bf Proof.} Since the letters $x_1,\ldots,x_n$ are super-letters,
every polynomial is a linear combination of monotonic super-words
by Lemma 7. Our present goal is to prove that the set of all
monotonic super-words is linearly independent. Let
\begin{equation}
\sum _i\alpha _iW_i=0
\label{ur}
\end{equation}

\noindent and assume that $W=[w_1]^{k_1}[w_2]^{k_2}\cdots
[w_m]^{k_m}$ is a leading super-word in (18). By Lemma 5, the
leading word of $W$ equals $w=w_1^{k_1}w_2^{k_2}\cdots w_m^{k_m}$.
Note that this word occurs exactly once in (18). Suppose, to the
contrary, that $W$ does also occur in the decomposition
$V=[v_1]^{m_1}[v_2]^{m_2}\cdots [v_k]^{m_k}$. Then the word $w$
is less than or equal to the leading word
$v=v_1^{m_1}v_2^{m_2}\cdots v_k^{m_k}$ in the decomposition of $V$,
which contradicts the fact that $W>V$ by Lemma 5.

Consider a free enveloping algebra in a given set of quantum
variables $H\langle x_1\ldots,x_n\rangle =G*{\bf k}\langle
x_1,\ldots,x_n\rangle $, on which the coproduct is defined by
setting
\begin{equation}
\Delta (x_i)=x_i\otimes 1+g_{x_i}\otimes x_i, \ \ \Delta (g) =
g\otimes g,
\label{izv}
\end{equation}

\noindent and group-like elements commute with variables via $xg=
\chi ^x(g)gx$. It follows that $G*{\bf k}\langle x_1,\ldots
,x_n\rangle$ turns into a Hopf algebra; for details, see [1, Sec.
3].

\smallskip
{\bf LEMMA 8.} {\it The coproduct at a super-letter $W=[w]$ is
represented thus:
\begin{equation}
\Delta ([w])=[w]\otimes 1+g_{w}\otimes [w]+\sum _i\alpha _i
G(W_i^{\prime \prime })W_i^{\prime }\otimes W_i^{\prime \prime },
\label{ko1}
\end{equation}

\noindent where $W_i^{\prime }$ are nonempty words in less
super-letters than is $[w]$ . Moreover, the sum of degrees of
super-words $W_i^{\prime }$ and $W_i^{\prime \prime }$ in each
variable $x_j$ equals the degree of $W$ in that variable, i.e., the
sum of structural elements of $W_i^{\prime}$ and $W_i^{\prime \prime
}$ is equal to the  constitution of} $W$.

{\bf Proof.} We use induction on the length of a word $w$. For
letters, by (19), there is nothing to prove. Let $W=[U,V]$,
$U=[u]$, and $V=[v]$. Assume that the decompositions
\begin{equation}
\Delta (U)=U\otimes 1+g_{u}\otimes U+\sum _i\alpha _i
G(U_i^{\prime \prime })U_i^{\prime }\otimes U_i^{\prime \prime },
\label{ko2}
\end{equation}

\noindent and
\begin{equation}
\Delta (V)=V\otimes 1+g_{v}\otimes V+\sum _j\beta _j
G(V_j^{\prime \prime })V_j^{\prime }\otimes V_j^{\prime \prime }
\label{ko3}
\end{equation}

\noindent satisfy the requirements of the lemma. Using (7) and
properties of a bicharacter $p$, we can write
$$ \Delta (W)=\Delta (U)\Delta (V)-p_{uv}\Delta (V)\Delta (U)= W\otimes 1+g_{w}\otimes W+ $$
$$ (1-p_{uv}p_{vu})g_uV\otimes U+\sum \beta _jp(U,V_j^{\prime \prime }) G(V_j^{\prime \prime 
})[UV_j^{\prime }]\otimes V_j^{\prime \prime }+ $$
\begin{equation}
\sum \beta _jg_uG(V_j^{\prime \prime })V_j^{\prime }\otimes (UV_j^{\prime \prime }-
p_{uv}p(V_j^{\prime },U)V_j^{\prime \prime }U)+
\end{equation}
$$\sum \alpha _iG(U_i^{\prime \prime})(U_i^{\prime }\cdot V-
p_{uv}p(V,U_i^{\prime \prime })V\cdot
U_i^{\prime })\otimes U_i^{\prime \prime }+
\sum \alpha _i p(U_i^{\prime },V)g_vG(U_i^{\prime \prime })U_i^{\prime }\otimes [U_i^{\prime \prime 
}V]+ $$
$$ \sum \alpha _i \beta _jG(U_i^{\prime \prime }V_j^{\prime \prime }) (p(U_i^{\prime }, V_j^{\prime 
\prime }) U_i^{\prime }V_j^{\prime }\otimes U_i^{\prime \prime }V_j^{\prime \prime }- 
p_{uv}p(V_j^{\prime },U_i^{\prime \prime }) V_j^{\prime }U_i^{\prime }\otimes V_j^{\prime \prime 
}U_i^{\prime \prime }). $$

Collecting similar terms in this formula will result in the
canceling of terms of the form $g_vU\otimes V$ only. We claim that
all left parts of the remaining tensors in (23) admit the required
decomposition. First, in view of the induction hypothesis, all
super-letters of all super-words $V_j^{\prime }$ are less than $V$,
which are in turn less than $W$ because $v$ is the end of a
standard word $w$. Moreover, by the induction hypothesis again,
$u$ cannot be the beginning of any word $u^{\prime }$ such that
the super-letter $[u^{\prime }]$ would occur in super-words
$U_i^{\prime }$. Therefore, $u>u^{\prime}$ implies $uv>u^{\prime
}$ or $W>[u^{\prime }]$. Thus all but the first and fourth
super-words on the left-hand sides of all tensors depend only on
super-letters which are less than $W$.

We want to apply Lemma 6 to the fourth tensor. Let $V_j^{\prime
}=\prod\limits _kV_{ik}$, where $V_{ik}=[v_{ik}]$ are less than
$V$. By formula (9) the polynomial $[U,V_j^{\prime }]$ is
a linear combination of words in the super-letters $V_{ik}$ and
skew commutators $[U,V_{ik}]$. By Lemma 6, each of these
commutators is a linear combination of words in the super-letters
$[v^{\prime }]$ such that $v^{\prime }\leq uv_{ik}$. In view of
$v_{ik}<v$, we obtain $v^{\prime }<uv=w$.

The statement concerning the constitution follows immediately from
formula (23) and the induction hypothesis.

\smallskip
{\bf LEMMA 9.} {\it The coproduct at a super-word $W$ is represented
thus:
\begin{equation}
\Delta (W)=W\otimes 1+G(W)\otimes W+\sum _i\alpha _i
G(W_i^{\prime \prime })W_i^{\prime }\otimes W_i^{\prime \prime },
\label{22}
\end{equation}

\noindent where the sum of constitutions  of $W_i^{\prime }$
and $W_i^{\prime \prime }$ equals the constitution of $W$.

{\bf Proof.} It suffices to observe that $\Delta $ is an
homomorphism of algebras. Here, we can no longer assert that}
$W_i^{\prime }<W$.

\medskip
\begin{center}
{\bf 3. BASIS FOR A CHARACTER HOPF ALGEBRA}
\end{center}
\smallskip

Consider a Hopf algebra $H$ generated by a set of skew primitive
semi-invariants $a_1,\ldots,a_n$ and by an Abelian group $G$ of
all group-like elements. Denote by $H_a$ a subalgebra generated by
$a_1,\ldots,a_n$. Then $H=GH_a$ since by definition,
semi-invariants obey the following commutation rule:
\begin{equation}
ag=\chi ^a(g)\cdot ga.
\label{kom}
\end{equation}

\noindent Let $x_1,\ldots,x_n$ be quantum variables with the same
parameters as $a_1,\ldots,a_n$, respectively, that is  $\chi
^{x_i}=\chi ^{a_i}$ and $g_{x_i}=g_{a_i}$. Then there exists an
homomorphism
\begin{equation}
\varphi :{\bf k}\langle x_1,\ldots,x_n\rangle \rightarrow H_a
\label{fi}
\end{equation}

\noindent which maps $x_i$ to $a_i$. This allows us to extend all
the combinatorial notions applied to the words in $x_1,\ldots,x_n$
in the above sections to the words in $a_1,\ldots,a_n$.

With $a_1,\ldots,a_n$ we associate the respective natural degrees
$d_1,\ldots,d_n$.$^*$\footnote{$^*$Note that further argument
will remain true for the case where $d_1,\ldots,d_n$ are arbitrary
positive elements of a linearly ordered additive Abelian group.}
 In this way, every word, super-letter, and super-word of 
a constitution $(m_1,\ldots ,m_n)$ have degree $m_1d_1+\ldots +m_nd_n$.

\smallskip
{\bf Definition 5.} A $G$-{\it super-word} is a product of the
form $gW$, where $g\in G$ and $W$ is a super-word. The degree,
constitution, length, and other concepts which apply with
$G$-super-words are defined by the super-word $W$. Alternatively, we
assume that the degree and the constitution of $g\in G$ are equal to
zero. In view of (25), every product of super-letters and
group-like elements equals a linear combination of $G$-super-words
of the same constitution.

\smallskip
{\bf Definition 6.} A super-letter $[u]$ is said to be {\it hard}
if it is not a linear combination of  words of the same degree in
less super-letters than is $[u]$ and of $G$-super-words of a lesser
degree.

\smallskip
{\bf Definition 7.} We say that the {\it height} of a super-letter
$[u]$ of degree $d$ equals a natural number $h$ if $h$ is least
with the following properties:

(1) $p_{uu}$ is a primitive root of unity of degree $t\geq 1$,
and either $h=t$ or $h=tl^k$, where $l$ is the characteristic
of the ground field;

(2) a super-word $[u]^h$ is a linear combination of super-words of
degree $hd$ in less super-letters than is $[u]$ and of
$G$-super-words of a lesser degree.

If, for the super-letter $[u]$, the number $h$ with the above
properties does not exist then we say that the height of $[u]$ is
infinite.

\smallskip
{\bf Definition 8.} The monotonic $G$-super-word
$$ g[u_1]^{n_1}[u_2]^{n_2}\cdots [u_k]^{n_k} $$

\noindent is said to be {\it restricted} if each of the numbers
$n_i$ is less than the height of the super-letter $u_i$.

\smallskip
{\bf THEOREM 2.} If a Hopf algebra $H$ is generated by a set
skew-primitive semi-invariants $a_1,\ldots, a_n$ and by an Abelian group
$G$ of all group-like elements, then the set of all monotonic
restricted $G$-super-words in hard super-letters constitute a
basis for $H$.

The proof will proceed through a number of lemmas. For
brevity, we call a super-word (a $G$-super-word) {\it admissible}
if it is monotonic restricted and is a word in hard super-letters
only.

\smallskip
{\bf LEMMA 10.} {\it Every nonadmissible super-word of degree $d$ is a
linear combination of lesser admissible super-words of degree $d$
and of admissible $G$-super-words of a lesser degree. Also, all
super-letters occurring in super-words of degree $d$ of this linear 
combination are less than or equal to a greatest super-letter 
of the super-word given}.

The {\bf proof} is by induction on the degree. Assume that the lemma
is valid for super-words of degree $<m$. Let $W$ be a least
super-word of degree $m$ for which the required representation
fails. By Lemma 7, the super-word $W$ is monotonic. If it has a
nonhard super-letter, by definition, we can replace it with a
linear combination of $G$-super-words of a lesser degree and of words
in less super-letters of the same degree. Removing the parentheses
turns $W$ into a linear combination of $G$-super-words of a lesser
degree and of lesser super-words of the same degree, a contradiction
with the choice of $W$. If $W$ contains a subword $[u]^k$,
where $k$ equals the height of $[u]$, then we can replace it as
is specified above, which gives us a contradiction again. Thus the
$W$ is itself monotonic restricted and is a word in hard
super-letters only.

In order to prove Theorem 2, it remains to show that admissible
$G$-super-words are linearly independent. Consider an arbitrary linear
combination $\bf T$ of admissible $G$-super-words and let
$U=U_1^{n_1}U_2^{n_2}\cdots U_k^{n_k}$ be its leading super-word of
degree $m$. Multiplying, if necessary, that combination by a
group-like element, we can assume that $U$ occurs once without a
group-like element:
\begin{equation}
{\bf T}=U+\sum _{j=1}^r\alpha _jg_jU+\sum _i\alpha _ig_iV_{i1}^{n_{i1}}
V_{i2}^{n_{i2}}\cdots V_{is}^{n_{is}}.
\label{lko1}
\end{equation}

In the next three lemmas, we accept the following
inductive assumption on $m$ and on $r$:
$$
(*) \ \ \ \matrix{\mbox{the set of all admissible} \
G\mbox{-super-words of degree} \ m \ \mbox{which are less than} \ U,
\hfill \cr
\mbox{of admissible} \
G\mbox{-super-words of degree} <m, \mbox{and of} \ G\mbox{-super-words}
\hfill \cr
g_jU, \ 1\leq j\leq r, \
\mbox{is linearly independent [we can assume} \
r=0 \mbox{ in (27)].} \hfill}
$$

\noindent In view of this assumption and Lemma 10, every super-word
of degree $m$ which is less than $U$, and every super-word of
degree $<m$, can be {\it uniquely} decomposed into a linear
combination of admissible $G$-super-words. For brevity, such will
be referred to as a {\it basis} decomposition.

\smallskip
{\bf LEMMA 11.} {\it If $\bf T$ is a skew primitive element then $r=0$
and all $G$-super-words of degree $m$ in (27) are super-words}.

{\bf Proof.} Rewrite the linear combination $\bf T$ as follows:
\begin{equation}
{\bf T}=U+\sum _{i=1}^k\gamma _i g_iW_i+W^{\prime },
\label{u1}
\end{equation}

\noindent where $g_iW_i$ are distinct $G$-super-words
of degree $m$ in (27) and $W^{\prime }$ a linear combination of
$G$-super-words of degree $<m$. In the expression
\begin{equation}
\Delta ({\bf
T})-{\bf T}\otimes 1-g_t\otimes {\bf T}, \label{pri}
\end{equation}

\noindent consider all tensors of the form $gW\otimes \ldots \ $,
where $W$ is of degree $m$. By Lemma 9, the sum of all such
tensors equals
\begin{equation}
\sum _{i=1}^r\gamma _ig_iW_i\otimes g_i-\sum _{i=1}^r
\gamma _ig_iW_i\otimes 1=
\sum _{i=1}^r \gamma _ig_iW_i\otimes (g_i-1).
\label{su}
\end{equation}

\noindent By inductive assumption $(*)$, the elements $g_iW_i$ are
linearly independent modulo all left parts of tensors of degree $<m$
in (29). Therefore, if (29) vanishes then either $\gamma _i=0$ or
$g_i=1$ for every $i$, as required.

\smallskip
{\bf LEMMA 12.} {\it If $\bf T$ is a skew primitive element then
$U=U_1^{n_1}$ and all super-words of degree $m$ except $U$ are words
in less super-letters than is} $U_1$.

{\bf Proof.} By the preceding lemma, we can assume that
\begin{equation}
{\bf T}=
\sum _i\alpha _ig_iV_{i1}^{n_{i1}}
V_{i2}^{n_{i2}}\cdots V_{is}^{n_{is}},
\label{lko}
\end{equation}

\noindent where $V_{ij}=[v_{ij}]$ are hard super-letters, $\alpha
_i$ are nonzero coefficients, and $g_i=1$ if $V_i$ is of degree
$m$. We apply coproduct to (31). By (8), then, the right-hand side
assumes the form
\begin{equation}
\sum _i\alpha _i (g_i\otimes g_i)\prod _{j=1}^s(V_{ij}\otimes 1+g_{ij} \otimes V_{ij}+
\sum _{\theta }g_{ij\theta }V_{ij\theta }^{\prime }\otimes
V_{ij\theta }^{\prime \prime })^{n_{ij}},
\label{lk1}
\end{equation}

\noindent where $V_{ij\theta }^{\prime }<V_{ij}$ and ${\rm
deg}V_{ij\theta }^{\prime }+{\rm deg}V_{ij\theta }^{\prime \prime
}={\rm deg}V_{ij}$.

Let $[v]$ be the largest super-letter occurring in super-words of
degree $m$ in (31). Since all super-words of (31) are monotonic,
this super-letter stands at the end of some super-words $V_i$,
i.e., $[v]=V_{is}$. If one of these super-words depends only on
$[v]$, i.e., $V_i=[v]^k$, then $V_i$ is a leading term, as
required. Therefore, we assume that every super-word of degree $m$
ending with $[v]$ is a word in more than one super-letter.

Let $k$ be the largest exponent $n_{is}$ of $[v]$ in ${\bf T}$.
Consider all tensors of the form $g[v]^k\otimes \ldots $ obtained
in (32) by removing the parentheses and applying the basis
decomposition to all left parts of tensors in all terms except ${\bf
T}\otimes 1$ (all of these terms are of degree $<m$).

All left parts of tensors which appear in
$$
\Delta (V_i)=(g_i\otimes g_i)\prod _{j=1}^s(V_{ij}\otimes 1+g_{ij}
\otimes V_{ij}+
\sum _{\theta }g_{ij\theta }V_{ij\theta }^{\prime }\otimes
V_{ij\theta }^{\prime \prime })^{n_{ij}}
$$

\noindent by removing the parentheses arise from the word
$V_i=\alpha _ig_iV_{i1}^{n_{i1}} V_{i2}^{n_{i2}}\cdots
V_{is}^{n_{is}}$ by replacing some of the super-letters $V_{ij}$
either with group-like elements $g_{ij}$ or with $G$-super-words
$g_{ij \theta }V_{ij\theta }^{\prime }$ of a lesser degree in less
super-letters. The right parts are, respectively, products obtained
by replacing super-letters $V_{ij}$ or super-words $V_{ij\theta
}^{\prime \prime}$ multiplied from the left by $g_i$.

If, under the replacements above, a new super-word is greater in
degree than $[v]^k$, then its basis decomposition
will give rise to terms of the form $g[v]^k\otimes \ldots \, $. In
this case, however, the right parts of those terms are of degree
less than $m-k{\rm deg}([v])$ since the sum of degrees of both parts
of the tensors either remains equal to $m$ or decreases.

If a new super-word is of degree less than the degree of $[v]^k$,
or the super-word is itself less than $[v]^k$ then its basis
decomposition will be freed of terms of the form $g[v]^k \otimes
\ldots $; see Lemma 10.

If a new super-word is of degree equal to that of $[v]^k$ and $V_i$
is of degree less than $m$ then the new super-word can be greater
than or equal to $[v]^k$. In this case the right-hand sides of the
new tensors are of degree less than $m-k{\rm deg}([v])$ because the
sum of degrees of the left- and right-hand sides of the tensors is
less than $m$.

If a new super-word is of degree equal to the degree of $[v]^k$,
but $V_i$ does not end with $[v]^k$, i.e., $V_i=W_i[v]^s$,
$0\leq s<k$, then the new super-word is less than $[v]^k$
since its first super-letter is less than $[v]$. (All
super-letters of $W_i$ cannot be replaced with group-like elements,
since otherwise the new word would be of degree less than or equal
to the degree of $[v]^s$.)

Finally, if $V_i=W_i[v]^k$ then a super-word of degree $k{\rm
deg}([v])$, which is greater than or equal to $[v]^k$, may
appear only if all super-letters of the super-words $W_i$ are
replaced with group-like elements, but $[v]$ is not. Here, the
resulting tensor is of the form $g(W_i)[v]^k\otimes \alpha _iW_i$.

We fix an index $t$ such that $V_t$ ends with $[v]^k$, letting
$t=1$. Then the sum of all tensors of the form $G(W_1)[v]^k\otimes
\ldots $ in $\Delta ({\bf T})-{\bf T}\otimes 1$ is equal to
\begin{equation}
G(W_1)[v]^k\otimes (\sum _j \alpha _jW_j+{\bf W}^{\prime}),
\label{ten}
\end{equation}

\noindent where ${\bf W}^{\prime}$ is a linear combination of basis
elements of degree less than $m-k{\rm deg}([v])$, and $j$ runs
through the set of all indices $i$ such that $V_i=W_i[v]^k$,
$G(W_i)=G(W_1)$, and the degree of $W_i$ equals $m-k{\rm
deg}([v])$. Since $W_i$ are distinct nonempty basis super-words
of degree less than $m$, tensor (33) is nonzero.

\smallskip
{\bf LEMMA 13.} {\it Under the conditions of Lemma $12,$ either $n_1=1$
or $p(U_1,U_1)$ is a primitive root of unity of degree $t\geq 
1$, in which case $n_1=t$ or the characteristic of a base
field equals $l>0$, and} $n_1=tl^k$.

{\bf Proof.} By the previous lemma, the linear combination ${\bf T}$
can be written in the form
\begin{equation}
{\bf T}=U^k+ \sum _i\alpha _ig_iV_{i1}^{n_{i1}}
V_{i2}^{n_{i2}}\cdots V_{is}^{n_{is}}, \label{lko11}
\end{equation}

\noindent where $U=[u]$ is greater than all super-letters
$V_{ij}$ for $V_i$ of degree $m$. First let $\xi
=1+p_{uu}+p_{uu}^2+ \ldots +p_{uu}^{k-1}\neq 0$ and assume $k>1$.

In the basis decomposition of $\Delta ({\bf T})-{\bf T}\otimes 1$,
consider tensors of the form $U^{k-1}\otimes \ldots \, $. All
super-letters $V_{ij}$ in super-words of degree $m$ are less than
$[u]$; therefore, tensors of this form may appear under the basis
decomposition of a tensor of $\Delta (V_i)-V_i\otimes 1$ only if
either the left part of that tensor is of degree greater than
$(k-1){\rm deg}([u])$ or $V_i$ is of degree less than $m$. In
either case the right part is of less degree than is $[u]$.
As above, if we remove the parentheses in
\begin{equation}
\Delta (U^k)=(U\otimes 1+g_u\otimes U+
\sum _{\tau }U_{\tau }^{\prime }\otimes U_{\tau }^{\prime \prime 
})^k ,\label{uk1} \end{equation}

\noindent we see that the left parts of the resulting tensors
arise from the super-word $U^k$ by replacing some super-letters $U$
either with $g_u$ or with super-words $U_{\tau }^{\prime }$ of a
lesser degree in less super-letters than is $U$. It follows that a
super-word of degree $(k-1){\rm deg}(U)$ which is greater than or
equal to $U^{k-1}$ appears only if exactly one super-letter is
replaced with a group element. Using the commutation rule
$U^sg_u=p_{uu}^sg_uU^s$, we see that the sum of all tensors of the
form $g_uU^{k-1}\otimes \ldots $ equals
$$ g_uU^{k-1}\otimes (\xi U+{\bf W}), $$

\noindent where ${\bf W}$ is a linear combination of basis
$G$-super-words of degree less than ${\rm deg}(U)$. Consequently,
(29) is nonzero for $k\neq 1$.

Now let $\xi =0$. Then $p_{uu}^k=1$. Therefore, $p_{uu}$ is a
primitive root of unity of some degree $t$ (we put $t=1$ if
$p=1$) and the number $k$ is divisible by $t$. We can write $k$
in one of the forms $t\cdot q$ or $tl^k\cdot q$, where $l$ is
the characteristic of a base field, in which $q\cdot 1\neq
0$. Put $h=t$ or $h=tl^k$, respectively. Since $(U\otimes
1)\cdot (g_u\otimes U)= p_{uu} (g_u\otimes U)\cdot (U\otimes 1)$,
use will be made of the quantum Newton binomial formula
\begin{equation}
(U\otimes 1+g_u\otimes U)^h=U^h\otimes
1+g_{uu}^h\otimes U^h. \label{m1}
\end{equation}

\noindent This implies that if we remove the parentheses in
\begin{equation}
\Delta (U^h)=((U\otimes 1+g_{uu}\otimes U)+
\sum _iG(U_i^{\prime \prime })
U_i^{\prime }\otimes U_i^{\prime \prime })^h,
\label{co11}
\end{equation}

\noindent then Lemma 8 gives
\begin{equation}
\Delta (U^h)=U^h\otimes 1+g_{uu}^h\otimes U^h+
\sum _{\theta }G(U_{\theta }^{\prime \prime })
U_{\theta }^{\prime }\otimes U_{\theta }^{\prime \prime },
\label{co2}
\end{equation}

\noindent where all super-words $U_{\theta }^{\prime }$ are less
than $U^h$ and are of less degree than is $U^h$. In this
formula, we note, all terms $U^r\otimes \ldots $, $r<h$, whose left
parts are greater than $U^h$, are banished. This allows us to treat
$U^h$ in (34) as a single block, or as a new formal super-letter
$\{ U^h\} $ such that $\{ U^h\}<U$, and $\{ U^h\}>[v_{ij}]$ if
$u^h>v_{ij}$ (which is equivalent to $u>v_{ij}$ by Lemma 3), i.e.,
\begin{equation}
{\bf T}=\{ U^h\} ^q+ \sum _i\alpha _ig_iV_{i1}^{n_{i1}}
V_{i2}^{n_{i2}}\cdots V_{is}^{n_{is}}.
\label{lko5}
\end{equation}

\noindent Since $p(U^h,U^h)=p^{h\cdot h}=1$, we have
$$ \xi _1=1+p(U^h,U^h)+\ldots +p(U^h,U^h)^{q-1}=q\neq 0. $$

\noindent As in the case above, assuming that $\{ U^h\} $ is a
single block, we can compute the sum of all tensors of the form
$g_u^h\{ U^h\} ^{q-1}\otimes \ldots $ in the basis decomposition of
$\Delta ({\bf T})-{\bf T}\otimes 1$ (provided that $q>1$):
\begin{equation}
g_{uu}^h\{ U^h\} ^{q-1}\otimes
(q\cdot \{ U^h\}+ {\bf W}),
\label{ten1}
\end{equation}

\noindent where ${\bf W}$ is a linear combination of basis
$G$-super-words of less degree than is $U^h$. By the
induction hypothesis, tensor (40) is nonzero, and so therefore is
(29).

The equality ${\bf T}=0$ does not hold. Indeed, if it did, then
${\bf T}$ would be a skew primitive element, which is nonzero in
view of Lemma 13 and definitions of hard super-letters and their
heights. Inductive assumption $(*)$ for $r=0$ is obviously valid
if $U$ is smallest among generators $a_i$, since group-like
elements, i.e., $G$-super-words of degree zero, are always
linearly independent. Theorem 2 is proved.

\medskip
\begin{center}
{\bf 4. SOME COROLLARIES}
\end{center}
\smallskip

In this section, again we write $H$ for a Hopf algebra generated by
an Abelian group $G$ of all group-like elements and by skew
primitive semi-invariants $a_1,\ldots,a_n$ with which degrees
$d_1,\ldots,d_n$ are associated.

\smallskip
{\bf COROLLARY 1.} {\it The set of all $G$-words in $a_1,\ldots,a_n$,
obtained by dropping all brackets from monotonic restricted
$G$-super-words in hard super-letters, constitute a basis for} $H$.

{\bf Proof.} Decompose an arbitrary word $v$ in $a_1,\ldots,a_n$
as is specified in Theorem 1, namely, $v=\sum\limits _j\alpha
_jV_j$, where $V_j=[v_{j1}]^{n_1}\cdots [v_{jk}]^{n_k}$ are
monotonic super-words of the same constitution. By Lemma 5, the leading
word appearing in $\sum\limits _j\alpha _jV_j$ under decomposition
(7) equals $v_s=v_{s1}^{n_1}\cdots v_{sk}^{n_k}$, where $V_s$ is
the leading super-word among all $V_j$. Therefore, $v=v_s$,
$\alpha _s=1$ --- this is still a decomposition in the free
algebra.

We use induction on the degree. Let $w$ be a minimal word of degree
$d$ which is not a linear combination of the $G$-words specified
in the statement. The word, as in the preceding paragraph, is
decomposed thus: $w=\sum\limits _j\alpha _jW_j$. If the leading
super-word $W_s$ is admissible, then $w$ arises from $W_s$ by
dropping the brackets, and so there is nothing to prove. If $W_s$
is not admissible, $W_s$ is the required linear combination by
Lemma~10 and inductive assumption $(*)$. We have $w=(\sum\limits
_{j\neq s}\alpha _jW_j)+W_s$, where the first summand is a linear
combination of words which are less than $w$, and again the
inductive assumption applies.

We argue for linear independence. Let $\sum\limits _{it}\beta
_{it}g_{it}w_i=0$. Then $w_i=\sum \alpha _{ij}W_{ij}$, where
$w_i$ is obtained by dropping all brackets from the leading
super-word $W_{is}$, and $\alpha _{is}=1$. Therefore, $W_{is}$
are admissible super-words. Now the equality $\sum\limits
_{ijt}\beta _{it}\alpha _{ij}g_{it}W_{ij}=0$ leads us to a
contradiction. Indeed, by Lemma~10, the nonleading super-words
$W_{ij}$ decrease under the basis decomposition, either in degree
or in ordering.

Sometimes we find it useful to apply the following criterion which
allows us to forget about skew commutators in computing hard
super-letters.

\smallskip
{\bf COROLLARY 2.} {\it A super-letter $[u]$ is hard if and only if the
standard word $u$ is not a linear combination of lesser words of
degree ${\rm deg}(u)$ and of $G$-words of a lesser degree}.

{\bf Proof.} Let $u=\sum\limits_i\alpha _iw_i+u_0$, where $w_i<u$
and ${\rm deg}(u_0)<{\rm deg}(u)$. Decompose the words $u$ and
$w_i$ as was done at the beginning of Corollary 1. We obtain
$u=[u]+\sum\limits _j\beta _jU_j$ and $w_i=\sum\limits _{t}\beta
_tW_{it}$, where the super-words $U_j$ are less than $[u]$, and
$w_i$ equals the leading word of a polynomial defined by the
leading super-word $W_{is}$. Since $u>w_i$, we have
$[u]>W_{is}$ by Lemma 5, and hence $[u]$ is greater than all
$W_{it}$. Consequently, the basis decompositions of $U_i$ and
$W_{it}$ have only super-words which either are less than $[u]$ or
of a lesser degree. For hard super-letters $[u]$, therefore, the
equality $[u]+\sum\limits _j\beta _jU_j- \sum\limits _i\alpha
_i\sum\limits _{t}\beta _tW_{it}-u_0=0$ is an impossibility.

Conversely, if $[u]=\sum \alpha _i W_i+U_0$, where $W_i$ depends
on super-letters less than $[u]$, then
$$ u=[u]+(u-[u])=\sum \alpha _i W_i+U_0+(u-[u]),$$

\noindent and the polynomial in the right part has no monomials
whose degree equals the degree of $u$ and which are greater than or
equal to $u$.

\smallskip
{\bf COROLLARY 3.} {\it Lemmas $11$-$1$3 are valid without assumption} $(*)$.

\smallskip
{\bf COROLLARY 4.} {\it A Hopf algebra $H$ is finite-dimensional if and
only if the group $G$ and the set of all hard super-letters are
finite, and each hard super-letter has finite height}.

\smallskip
{\bf COROLLARY 5.} {\it If $H$ has only a finite number of hard
super-letters and $G$ is finitely generated, then $H$ is (left and
right) Noetherian}.

\smallskip
{\bf COROLLARY 6.} {\it Let the group algebra ${\bf k}[G]$ has no
zero divisors. If $H$ has only a finite number of hard
super-letters, of which each has infinite height, then $H$ has
no zero divisors and has a classical skew field of quotients}.

As in the case of classical Lie algebras, in order to prove the last
two corollaries, we need only construct on $H$ a filtration
$H_0\subseteq H_1\subseteq \ldots \subseteq H_k\subseteq \ldots$
such that the associated graded algebra $D(H)$ satisfies the
required properties; see, e.g., [17, Ch. V, Sec. 3, Thms. 4, 5].

\smallskip
{\bf Construction of the filtration.} Assume that $H$ has finitely
many hard super-letters. Consider a set ${\bf R}$ of all words in
$a_1,\ldots, a_n$, whose degree does not exceed the maximal
degree of a hard super-letter multiplied by a maximal finite height,
or by 2 if all heights are infinite. In this case ${\bf R}$ is
composed of all standard words defining hard super-letters and of
all products $uv$, $u^h$, where $[u]$ and $[v]$ are hard
super-letters and $h$ is the height of $u$. Let words of ${\bf
R}$ all respect the lexicographical ordering described at the
beginning of Sec. 1 and $n(u)$ be the number of words in ${\bf R}$
which are less than or equal to $u$. The largest word $a_1$ is
defined by the number $L=n(a_1)$. Denote by $M$ an arbitrary
natural number which is greater than the length of any word in
${\bf R}$. Define the filtration degree on hard super-letters
using the formula
\begin{equation}
{\rm Deg}([u])=M^{L+1}{\rm deg}(u)+M^{n(u)},
\label{De}
\end{equation}

\noindent where ${\rm deg}(u)$ is specified by the constitution of $u$,
and ${\rm deg}(u)=d_1m_1+\ldots+d_nm_n$. The filtration degree of a
basis element $gW$ equals the sum of filtration degrees of all of
its super-letters. The filtration degree of an arbitrary element
${\bf T}\in H$ equals the maximal filtration degree of the basis
elements occurring in its basis decomposition.

\smallskip
{\bf LEMMA 14.} {\it The function ${\rm Deg}$ defines a filtration on
$H$, so that}
$$ H_0={\bf k}[G],\quad H_k=\{ {\bf T}\in H\, |\, {\rm Deg}({\bf T})\leq k\}. $$

{\bf Proof.} We have to show that $H_kH_s\subseteq H_{k+s}$, i.e.,
${\rm Deg}({\bf T}_1\cdot {\bf T}_2)\leq {\rm Deg}({\bf
T}_1)+{\rm Deg}({\bf T}_2)$. To do this, we construct an
additional degree function $D^{\prime }$ on a set of all linear
combinations of (not necessarily admissible) super-words in the
super-letters defined by all standard words of the vocabulary ${\bf
R}$ The $D^{\prime }$-degree of a super-letter is defined by
formula (41). The $D^{\prime }$-degree of a product of
super-letters equals the sum of degrees of its factors.
Accordingly, the $D^{\prime }$-degree of a linear combination
equals the maximum of $D^{\prime }$-degrees of its summands. Of
course we do not claim that the various linear combinations defining
equal elements of $H$ have the same $D^{\prime }$-degrees.

If we assume that ${\bf T}_1 \cdot {\bf T}_2$ is obtained from
${\bf T}_1$ and ${\bf T}_2$ by merely removing the parentheses,
then ${\rm Deg}({\bf T}_1)+{\rm Deg}({\bf T}_2)= D^{\prime }({\bf
T}_1\cdot {\bf T}_2)$. Therefore, it suffices to specify how a
basis decomposition of super-words proceeds in a way that $D^{\prime
}$ is kept unincreased. Our plan is as follows. First, we replace
nonhard super-letters via Definition 6, next replace all subwords
$[u]^h$, where $h$ is the height of a hard super-letter $[u]$,
then apply the decreasing decomposition of Lemma 7, and again
replace nonhard super-letters, etc.

Let $[u]$ be a nonhard super-letter defined by $u$ in ${\bf R}$
as follows:
\begin{equation}
[u]=\sum _i\alpha _i \prod _j[w_{ij}]+\sum _s\alpha _sg_s\prod _t[v_{st}],
\label{um1}
\end{equation}

\noindent where $[w_{ij}]$ are less than $[u]$,
$n(w_{ij})\leq n(u)-1$, of the combined ordinary degree ${\rm
deg}(u)$, and $\sum\limits _t{\rm deg}(v_{st})\leq {\rm
deg}(u)-1$. By the definition of $R$, all words $w_{ij}$ and
$v_{st}$, as well as the words which result from the summands of
(42) by dropping all brackets and multiplication signs, belong to
the vocabulary ${\bf R}$. In particular, their lengths are less
than $M$. Therefore, the number of factors in every summand of
(42) is less than $M$. Thus we may write
\begin{equation}
D^{\prime }(\prod _j[w_{ij}])=\sum _jM^{L+1}{\rm deg}(w_{ij}) +\sum _j M^{n(w_{ij})}<
M^{L+1}{\rm deg}(u)+M\cdot M^{n(u)-1}=D^{\prime }(u).
\end{equation}

\noindent Similarly,
\begin{equation}
\begin{array}{c}
 D^{\prime }(\prod _t[v_{st}])=\sum _tM^{L+1}{\rm deg}(v_{st}) +
\sum _t M^{n(v_{st})}< \\
M^{L+1}({\rm deg}(u)-1)+M\cdot M^L=
D^{\prime }(u)-M^{n(u)}-M^{L+1}+M^{L+1}<D^{\prime }(u).
\end{array}
\label{ume1}
\end{equation}

The argument is the same if $[u]^h$ is changed with $\sum\limits_i
\prod\limits_j[w_{ij}]+\sum\limits_sg_s\prod\limits_t[v_{st}]$,
where $[u]$ is a hard super-letter of height $h$. We have
\begin{equation}
D^{\prime }(\sum _i \prod _j[w_{ij}]+\sum _sg_s\prod _t[v_{st}])<
D^{\prime }([u]^h).
\label{10}
\end{equation}

We proceed to the decreasing process of Lemma 7. The second summand
in (17) has the same $D^{\prime }$-degree as $W$. By Lemma 6,
the factor $[UU_1]$ in the first summand can be represented as
$\sum\limits_i\alpha_i\prod\limits_j[w_{ij}]$, where the
super-letters $[w_{ij}]$ are less than $U$, and they all are in
${\bf R}$ since their ordinary degrees are less than or equal to
the ordinary degree of $UU_1$. Therefore, $n(w_{ij})\leq 
n(u)-1$. And the constitutions being equal indicates that the number
of factors in $\prod\limits_j$ does not exceed $M$. We have
\begin{equation}\begin{array}{c}
 D^{\prime }(\prod _j[w_{ij}])= \sum _j M^{L+1}{\rm deg}(w_{ij})+
\sum M^{n(w_{ij})}< \\
M^{L+1}{\rm deg}(UU_1)+M\cdot M^{n(u)-1}<
M^{L+1}{\rm deg}(UU_1)+M^{n(u)}+M^{n(u_1)}=D^{\prime }(UU_1).
\end{array}
\end{equation}

\noindent The lemma is proved.

Note that the strict inequality signs in (46) show that the
inequality
\begin{equation}
{\rm Deg}([u]\cdot [v]-p_{uv}[v]\cdot [u])<{\rm Deg}(u)+{\rm Deg}(v)
\label{12}
\end{equation}

\noindent is valid for all hard super-letters $U=[u]$, $U_1=[v]$,
$u>v$. Similarly, (45) gives us
\begin{equation}
{\rm Deg}([u]^h)<h{\rm Deg}(u),
\label{14}
\end{equation}

\noindent where $h$ is the height of a hard super-letter $[u]$.

\smallskip
{\bf Associated graded algebra.} With each hard super-letter $[u]$
we associate a new variable $x_{u}$. Denote by $H^*$ an algebra
generated by $x_u$ and defined by the relations
$x_ux_v=p_{uv}x_vx_u$, where $u>v$. This algebra can be
constructed by iterating the skew polynomials construction. Namely,
let $[u_1]<[u_2]<\ldots <[u_s]$ all be hard super-letters. Denote
by $H^*_k$ a subalgebra generated by $x_{u_1},\ldots,x_{u_k}$.
Then $H^*_1$ is isomorphic to an algebra of polynomials in one
variable. The map $\varphi : x_v\rightarrow p_{uv}x_v$, where
$u=u_{k+1}$, $v=u_i$, $i\leq k$, determines an
automorphism of $H_k^*$. The commutation rule
$x_ux_v=p_{uv}x_vx_u$ can be written in the form
$x_ux_v=x_v^{\varphi }x_u$. Therefore, the algebra $H^*_{k+1}$
is isomorphic to an algebra of skew polynomials over $H^*_k$. In
particular, the algebra $H^*$ is Noetherian and has no zero
divisors.

Define the action of $G$ on $H^*$ by the formula $x_u^g=\chi
^u(g)\cdot x_u$. Let $H^*[G]$ be a skew group algebra of $G$
with coefficients from~$H^*$.

\smallskip
{\bf THEOREM 3.} {\it The associated graded algebra $D(H)$ is isomorphic
to the quotient algebra of $H^*[G]$ with respect to relations
$x_u^h=0$, where $[u]$ runs through the set of all hard
super-letters of finite height and $h$ equals the height of~$[u]$}.

{\bf Proof.} Denote by $x_u$ an element of $D(H)$ defined by the
coset $[u]+H_{m-1}$, where $[u]$ is a hard super-letter and
$m={\rm Deg}(u)$. Then the zero component ${\bf k}[G]$ and the
elements $x_u$ generate $D(H)$. Formulas (47) and (48) show that
$x_ux_v=p_{uv}x_vx_u$ and $x_u^h=0$ hold if $u>v$ and $h$ is the
height of $[u]$. All the monotonic restricted $G$-words in
$x_u$ are linearly independent in $D(H)$ since the filtration
degree of every admissible word equals the sum of filtration
degrees of all super-letters of that word.

{\bf Braided bigraded Hopf algebras}. The above results can be
easily extended to the case of $(G, \lambda )$-graded Hopf
algebras (see, e.g., [18, p. 206]) and to the case of braided
bigraded Hopf algebras (cf. [1]). These objects are not Hopf
algebras in the ordinary sense, but still they are Hopf algebras in some
categories. In view of Radford's results, these algebras have
embeddings in the ordinary Hopf algebras; see [19, Thm. 1 and 
Prop.~2]. The embeddings are obtained by adding the elements of $G$ 
treated as group-like elements. In this case primitive elements
correspond to skew primitive ones. If a given $(G,\lambda )$-graded
Hopf algebra $H$ is generated by primitive elements then the
enveloping Hopf algebra is character and {\it separated}, i.e.,
$H_a\cap {\bf k}[G]=0$.

All the concepts of this article can be easily extended to the case
of $(G,\lambda )$-graded Hopf algebras by vanishing group-like
elements and by replacing the bicharacter $p(u,v)$ with $\lambda
^{-1}(g_v,g_u)$. We are now in a position to formulate the
following corollarries.

\smallskip
{\bf COROLLARY 7.} {\it If a $(G,\lambda )$-graded Hopf algebra $H$ is
generated by a set of primitive elements $a_1,\ldots, a_n$ then the
set of all monotonic restricted words in hard super-letters
constitutes a basis for~$H$}.

\smallskip
{\bf COROLLARY 8.} {\it Let $H$ be a $(G,\lambda )$-graded Hopf
algebra generated by primitive elements. $H$ is
finite-dimensional if and only if the set of all hard super-letters
is finite and each hard super-letter has finite height}.

\smallskip
{\bf COROLLARY 9.} {\it Let $H$ be a $(G,\lambda )$-graded Hopf
algebra generated by primitive elements. If $H$ has only a finite
number of hard super-letters then it is Noetherian}.

\smallskip
{\bf COROLLARY 10.} {\it Let $H$ be a $(G,\lambda )$-graded Hopf
algebra which is generated by primitive elements and has only
a finite number of hard super-letters. If all these hard
super-letters have infinite height then $H$ has no zero divisors and
has a classical field of quotients}.

In the corollaries above, we note, restrictions on a group are
unnecessary since the group merely defines a bigrading, but the
algebra in question would not contain~it.

{\bf Acknowledgement}. I want to thank Dr. Juan Antonio Montaraz, Director
of FES-C, Dr. Suemi Rodriguez-Romo, and Virginia Lara Sagahon for
providing beautiful facilities for my research work at FES-C UNAM,
M\'exico. Thanks also are due to participants of Shirshov Seminar on Ring Theory
(Institute of Mathematics of RAS) for interesting comments on the subject matter.

\medskip
\begin{center}
{\bf REFERENCES}
\end{center}
\begin{enumerate}

\item V. K. Kharchenko, ``An algebra of skew primitive elements,''
{\it Algebra and Logic}, {\bf 37}, No. 2, 101--126 (1998).

\item Y. Ju. Reshetikhin, L. A. Takhtadzhyan, and L. D. Faddeev,
``Quantizations of Lie groups and Lie algebras,'' {\it Leningrad
Mat. Zh.}, {\bf 1}, No. 1, 193-225 (1990).

\item A. I. Shirshov, ``On free Lie rings,'' {\it Mat. Sb.}, {\bf
45}(87), No. 2, 113-122 (1958).

\item M. Hall, ``A basis for free Lie rings and higher commutators
in free groups,'' {\it Proc. Am. Math. Soc.}, {\bf 1}, 575-581
(1950).

\item A. I. Shirshov, ``On bases for free Lie algebra,'' {\it
Algebra Logika}, {\bf 1}, No. 1, 14-19 (1962).

\item A. I. Shirshov, ``Some algorithmic problems for Lie
algebras,'' {\it Sib. Mat. Zh.}, {\bf 3}, No. 2, 292-296 (1962).

\item K. T. Chen, R. H. Fox, and R. C. Lyndon, ``Free differential
calculus IV, the quotient groups of the lower central series,'' {\it
Ann. Math.}, {\bf 68}, 81-95 (1958).

\item A. I. Mal'tsev, ``On representations of models,'' {\it Dokl.
Akad. Nauk SSSR}, {\bf 108}, No. 1, 27-29 (1956).

\item M. Rosso, ``An analogue of the Poincare--Birkhoff--Witt
theorem and the universal $R$-matrix of $U_q(sl(N+1))$,'' {\it Comm.
Math. Phys.}, {\bf 124}, 307-318 (1989).

\item H. Yamane, ``A Poincar\`e-Birkhoff-Witt theorem for quantized
universal enveloping algebras of type $A_N$,'' {\it Publ., RIMS.
Kyoto Univ.}, {\bf 25}, 503-520 (1989).

\item G. Lusztig, ``Canonical bases arising from quantized
enveloping algebras,'' {\it J. Am. Math. Soc.}, {\bf 3}, No. 2,
447-498 (1990).

\item G. Lusztig, ``Quivers, perverse sheaves, and quantized
enveloping algebras,'' {\it J. Am. Math. Soc.}, {\bf 4}, No. 2,
365-421 (1991).

\item M. Kashiwara, ``On crystal bases of the q-analog of universal
enveloping algebras,'' {\it Duke Math. J.}, {\bf 63}, No. 2, 465-516
(1991).

\item E. J. Taft and R. L. Wilson, ``On antipodes in pointed Hopf
algebras,'' {\it J. Alg.}, {\bf 29}, 27-32 (1974).

\item J. W. Milnor and J. C. Moore, ``On the structure of Hopf
algebras,'' {\it Ann. Math. (2)}, {\bf 81}, 211-264 (1965).

\item A. I. Shirshov, ``Some problems in the theory of rings which
are close to associative,'' {\it Usp. Mat. Nauk}, {\bf 13}, No.
6(84), 3-20 (1958).

\item N. Jacobson, {\it Lie Algebras}, Interscience, New York
(1962).

\item S. Montgomery, {\it Hopf Algebras and Their Actions on Rings,
Reg. Conf. Ser. Math.}, {\bf 82}, Am. Math. Soc., Providence, RI
(1993).

\item D. E. Radford, ``The structure of Hopf algebras with a
projection,'' {\it J. Alg.}, {\bf 92}, 322-347 (1985).

\end{enumerate}
\end{document}